\documentclass[11pt]{amsart}
\def\Q{{\mathbf Q}}
\def\Z{{\mathbf Z}}

\def\F{{\mathbf F}}
\def\SS{{\mathbf S}}

\def\Sn{{\mathbf S}_n}
\def\An{{\mathbf A}_n}

\def\Gal{\mathrm{Gal}}

\def\End{\mathrm{End}}
\def\Aut{\mathrm{Aut}}

\def\fchar{\mathrm{char}}

\def\GL{\mathrm{GL}}
\def\SL{\mathrm{SL}}

\def\dim{\mathrm{dim}}

\newtheorem{thm}{Theorem}[section]
\newtheorem{lem}[thm]{Lemma}

\theoremstyle{definition}

\newtheorem{exs}[thm]{Examples}

\newtheorem{rem}[thm]{Remark}

\title[Hyperelliptic jacobians]
{Hyperelliptic jacobians without complex multiplication in positive
characteristic}
\author[Yuri G. Zarhin]{Yuri G. Zarhin}
\address{Department of Mathematics, Pennsylvania State University,
University Park, PA 16802, USA}
\email{zarhin\char`\@math.psu.edu}
\thanks{Partially supported by the NSF}

\begin{document}
\maketitle
\section{Introduction}
The aim of this note is to prove that
in positive characteristic $p \ne 2$ the
jacobian $J(C)=J(C_f)$  of a hyperelliptic curve
$$C=C_f:y^2=f(x)$$
has only trivial endomorphisms
over an algebraic closure of the ground field $K$
if the Galois group $\Gal(f)$ of the polynomial
$f \in K[x]$ of even degree is ``very big".

More precisely, if $f$ is a polynomial of {\sl even} degree $n \ge 10$
and $\Gal(f)$ is either the symmetric group $\Sn$ or the
alternating group $\An$ then $\End(J(C))=\Z$.
Notice that it is known \cite{ZarhinMRL} that in this case (and even for all
integers $n\ge 5$) either $\End(J(C))=\Z$ or $J(C)$ is a supersingular
abelian variety and the real problem is how to prove that $J(C)$ is {\sl
not} supersingular.

There are some results of this type in the literature. Previously
Mori \cite{Mori1}, \cite{Mori2} has constructed explicit examples of
hyperelliptic jacobians without nontrivial endomorphisms. Namely, he proved
that if $K=k(z)$ is a field of rational functions in variable $z$ with
constant field $k$ of characteristic $p \ne 2$ then for each integer $g \ge
2$ the jacobian of a hyperelliptic $K$-curve
$$y^2=x^{2g+1}-x+z$$
has no nontrivial endomorphisms if $p$ does not divide $g(2g+1)$.

Our result stated above implies, in particular,  that for each integer $g
\ge 4$ the jacobian of a hyperelliptic $K$-curve
$$y^2=x^{2g+2}-x+z$$
has no nontrivial endomorphisms if $p$ does not divide $(g+1)(2g+1)$, since in this case the Galois group of
$x^{2g+2}-x+z$ over $K$ is $\SS_{2g+2}$.

\section{Main result}
\label{mainr}
Throughout this paper we assume that $K$ is a field of prime characteristic
$p$
different from $2$. We fix its algebraic closure $K_a$
and write $\Gal(K)$ for the absolute Galois group $\Aut(K_a/K)$.

\begin{thm}
\label{main}
Let $K$ be a field with $p=\fchar(K) > 2$,
 $K_a$ its algebraic closure,
$f(x) \in K[x]$ an irreducible separable polynomial of even degree $n
\ge 10$ such that the Galois group of $f$ is either $\Sn$ or $\An$.
Let $C_f$ be the hyperelliptic curve $y^2=f(x)$. Let  $J(C_f)$ be
its jacobian, $\End(J(C_f))$ the ring of $K_a$-endomorphisms of
$J(C_f)$. Then $\End(J(C_f))=\Z$.
\end{thm}

\begin{exs}
Let $k$ be a field of odd characteristic $p$. Let
$k(z)$ be the field of rational functions in variable $z$ with constant
field $k$. 

\begin{enumerate}
\item
Suppose $K_n=k(z_1, \cdots , z_n)$ is the field of rational functions in $n$ independent variables $z_1, \cdots , z_n$ over $k$. Then the Galois group of a polynomial $x^n-z_1 x^{n-1} + \cdots +(-1)^n z_n$ over $K_n$ is $\Sn$. Therefore if $n \ge 10$ is even then
the jacobian of the curve $y^2=x^n-z_1 x^{n-1} + \cdots +(-1)^n z_n$ has no nontrivial endomorphisms over an algebraic closure of $K_n$.
\item
Suppose
$h(x) \in k[x]$ is a {\sl Morse polynomial} of degree $n$  and $p$ does not
divide $n$. This means that the derivative $h'(x)$ of $h(x)$ has $n-1$ distinct roots $\beta_1, \cdots \beta_{n-1}$ (in an algebraic closure of $k$) and $h(\beta_i) \ne h(\beta_j)$ while $i\ne j$. For example, $h(x)=x^n-x$ enjoys these properties if and only if
$p$ does not divide $n(n-1)$.

 Then the Galois group of $h(x)-z$ over $k(z)$ is the symmetric
group $\Sn$ (\cite{Serre}, Th. 4.4.5, p. 41). Hence if $n\ge 10$ is even and $p$
does not divide $n(n-1)$ then
the jacobian of the curve $y^2=h(x)-z$ has no nontrivial endomorphisms over
an algebraic closure of $k(z)$. 
\item
Suppose $k$ is algebraically closed. 
Suppose $n=q+t$ where $q$ is a power of $p$ and $t>q$ is a positive integer
not divisible by $p$.
Then the Galois group of $x^n-xz^t+1$ over $k(z)$ is the alternating group
$\An$ (\cite{AbhMA}, Th.  1, p. 67). Clearly, if $t$ is odd then $n=q+t$ is
even and $n>2q\ge 6$, i.e., $n \ge 8$. In addition, $n\ge 10$ unless $q=3, t=5$. This implies that if $t$ is odd and $(q,t)\ne (3,5)$ then the jacobian of the curve $y^2=x^n-xz^t+1$ has no nontrivial
endomorphisms over an algebraic closure of $k(z)$.
\end{enumerate}
\end{exs}

As was already pointed out, in light of Th. 2.1 of \cite{ZarhinMRL},
 Theorem \ref{main} is an immediate corollary of the following auxiliary
statement.

\begin{thm}
\label{main2} Suppose   $n=2g+2$ is an even integer which is
greater than or equal to $10$. Suppose $f(x) \in K[x]$ is a
separable polynomial of degree $n$, whose Galois group
is either
$\An$ or $\Sn$. Suppose $C$ is the hyperelliptic curve $y^2=f(x)$ of
genus $g$ over $K$ and
 $J(C)$  is the jacobian of $C$.

Then $J(C)$ is not a supersingular abelian variety.
\end{thm}

\begin{rem}
\label{redA}
Replacing $K$ by its proper quadratic extension, we may  assume in the
course of the proof of Theorem \ref{main2} that $\Gal(f)=\An$. Also,
replacing $K$ by its abelian extension obtained by adjoining to $K$ all
$2$-power roots of unity, we may assume that $K$ contains all $2$-power
roots of unity.
\end{rem}

We prove Theorem \ref{main2} in the next Section.

\section{Proof of Theorem \ref{main2}}
So, we assume that $K$ contains all $2$-power roots of unity, $f(x) \in
K[x]$ is an irreducible separable polynomial of even degree $n=2g+2\ge 10$
and $\Gal(f)=\An$.
Therefore $J(C)$ is a $g$-dimensional
abelian variety defined over $K$.
The group $J(C)_2$ of its points of order $2$ is  a
$2g$-dimensional $\F_2$-vector space provided with the natural action of
$\Gal(K)$. It is well-known (see for instance \cite{ZarhinTexel})that the
image
of $\Gal(K)$ in $\Aut(J(C)_2)$ is canonically isomorphic to $\Gal(f)$.

Now Theorem  \ref{main2} becomes an immediate corollary of the following two
assertion.

\begin{lem}
\label{supernot}
Let $F$ be a field, whose characteristic is not $2$ and assume that $F$
contains all $2$-power roots of unity. Let $g$ be a positive integer and $G$
a finite simple non-abelian group enjoying the following properties:
\begin{enumerate}
\item[(a)]
Each nontrivial representation of $G$ in characteristic $0$ has dimension
$>2g$;
\item[(b)]
If $G' \to G$ is a surjective group homomorphism,
whose kernel is a central subgroup of order $2$ then each faithful
irreducible representation of $G'$ in characteristic zero has dimension
$ \ne 2g$.


\item[(c)]
Each nontrivial representation of $G$ in characteristic $2$ has dimension
$\ge 2g$.
\end{enumerate}

If $X$ is a $g$-dimensional abelian variety   over  $F$ such that the image
of $\Gal(F)$ in $\Aut(X_2)$ is isomorphic to $G$ then $X$ is not
supersingular.
\end{lem}

\begin{lem}
\label{repAn}
Suppose $n=2g+2\ge 10$ is an even integer.
Let us put $G=A_n$. Then
\begin{enumerate}
\item[(a)]
Each nontrivial representation of $G$ in characteristic $0$ has dimension
$\ge n-1>2g$;
\item[(b)]
Each nontrivial proper projective representation of $G$ in characteristic
$0$  has dimension $\ne 2g$;
\item[(c)]
Each nontrivial representation of $G$ in characteristic $2$ has dimension
$\ge 2g$.
\end{enumerate}
\end{lem}

Lemma \ref{supernot} will be proven in the next Section. We prove Lemma \ref{repAn} in Section \ref{reps}.

\section{Not supersingularity}
\label{main2p}
We keep all the notations of Lemma \ref{supernot}.
Assume that $X$ is supersingular. Our goal is to get a contradiction.
We write $T_2(X)$ for the $2$-adic Tate module of $X$
and
$$\rho_{2,X}:\Gal(F) \to \Aut_{\Z_2}(T_2(X))$$
for 
the corresponding $2$-adic representation. It is well-known that $T_2(X)$ is
a free $\Z_2$-module of rank $2\dim(X)=2g$ and
$$X_2=T_2(X)/2 T_2(X)$$
(the equality of Galois modules). Let us put
$$H=\rho_{2,X}(\Gal(F)) \subset \Aut_{\Z_2}(T_2(X)).$$
Clearly, the natural homomorphism
$$\bar{\rho}_{2,X}:\Gal(F) \to \Aut(X_2)$$
defining the Galois action on the points of order $2$ is the composition of
$\rho_{2,X}$ and (surjective) reduction map modulo $2$
$$\Aut_{\Z_2}(T_2(X)) \to \Aut(X_2).$$
This gives us a natural (continuous) {\sl surjection}
$$\pi:H \to \bar{\rho}_{2,X}(\Gal(F)) \cong G,$$
 whose kernel consists of elements of $1+2 \End_{\Z_2}(T_2(X))$. It follows
from the property \ref{supernot}(c) and equality $\dim_{\F_2}(X_2)=2g$ that
the $G$-module $X_2$ is absolutely simple and therefore the $H$-module $X_2$
is also absolutely simple. Here the structure of $H$-module is defined on
$X_2$ via
$$H\subset\Aut_{\Z_2}(T_2(X)) \to \Aut(X_2).$$
Let $V_2(X)=T_2(X)\otimes_{\Z_2}\Q_2$ be the $\Q_2$-Tate module of $X$. It is
well-known that $V_2(X)$ is the $2g$-dimensional $\Q_2$-vector space and
$T_2(X)$ is a $\Z_2$-lattice in $V_2(X)$.
This implies easily that the $\Q_2[H]$-module $V_2(X)$ is also absolutely
simple.

The choice of polarization on $X$ gives rise to a non-degenerate alternating
bilinear form (Riemann form) \cite{MumfordAV}
$$e: V_{2}(X) \times V_2(X) \to \Q_2(1) \cong \Q_2.$$
Since $F$ contains all $2$-power roots of unity, $e$ is $\Gal(F)$-invariant
and therefore is $H$-invariant. In particular,
$$H \subset \SL(V_2(X)).$$

There exists a finite Galois extension $L$ of $K$ such that all
endomorphisms of $X$ are defined over $L$. We write $\End^0(X)$ for the
$\Q$-algebra $\End(X)\otimes\Q$ of endomorphisms of $X$.
Since $X$ is supersingular,
$$\dim_{\Q}\End^0(X)=(2\dim(X))^2=(2g)^2.$$
Recall (\cite{MumfordAV}) that the natural map
$$\End^0(X)\otimes_{\Q}\Q_{2} \to \End_{\Q_{2}}V_{2}(X)$$
is an embedding.
Dimension arguments imply that
$$\End^0(X)\otimes_{\Q}\Q_{2} = \End_{\Q_{2}}V_{2}(X).$$
Since all endomorphisms of $X$ are defined over $L$, the image
$$\rho_{2,X}(\Gal(L)) \subset \rho_{2,X}(\Gal(F)) \subset\Aut_{\Z_2}(T_2(X))
\subset\Aut_{\Q_2}(V_2(X)$$
commutes with $\End^0(X)$. This implies that
$\rho_{2,X}(\Gal(L))$ commutes with  $\End_{\Q_{2}}V_{2}(X)$ and therefore
consists of scalars. Since
$$\rho_{2,X}(\Gal(L)) \subset \rho_{2,X}(\Gal(F)) \subset \SL(V_2(X)),$$
$\rho_{2,X}(\Gal(L))$ is a finite group. Since $\Gal(L)$ is a subgroup of
finite index in $\Gal(F)$, the group $H=\rho_{2,X}(\Gal(F))$ is also finite.
In particular, the kernel of the reduction map modulo $2$
$$\Aut_{\Z_2}T_2(X) \supset H \to G \subset \Aut(X_2)$$
consists of periodic elements and, thanks to Minkowski-Serre Lemma,
$Z:=\ker(H \to G)$ has exponent $1$ or $2$. In particular, $Z$ is
commutative. Since
$$Z \subset H \subset \SL(V_2(X)),$$
$Z$ is a $\F_2$-vector space of dimension $d<2g$.
This implies that the adjoint action
$$G \to \Aut(Z) \cong\GL_d(\F_2)$$
is trivial, in light of property \ref{supernot}(c). This means that $Z$ lies
in the center of $H$. Since the $\Q_2[H]$-module $V_2(X)$ is faithful absolutely
simple, $Z$ consists of scalars. This implies that either $Z=\{1\}$ or
$Z=\{\pm 1\}$. If $Z=\{1\}$ then $H \cong G$ and $V_2(X)$ is a faithful
$\Q_2[G]$-module of dimension $2g$ which contradicts the property
\ref{supernot}(a). Therefore $Z=\{\pm 1\}$ and $H \to G$
is a surjective group homomorphism,
whose kernel is a central subgroup of order $2$.
But $V_2(X)$ is a faithful absolutely simple $\Q_2[H]$-module of dimension
$2g$ which contradicts the property \ref{supernot}(b).
This ends the proof of Lemma \ref{supernot}.

\section{Representation theory}
\label{reps}
\begin{proof}[Proof of Lemma \ref{repAn}]
The property (a) follows easily from Th. 2.5.15 on p. 71 of \cite{JK}. The
property (c) follows readily from Th. 1.1 on p. 127 of \cite{Wagner}.

In order to prove the property (b), recall (Th. 1.3(ii) on pp. 583--584 of
\cite{Wagner2}) that  each proper projective representation of $\An$ in characteristic $\ne 2$ has dimension divisible
by $2^{[(n-s-1)/2]}$ where $s$ is the exact number of terms in the dyadic expansion of $n$. We have $n=2^{w_1}+ \cdots + 2^{w_s}$ where  $w_i$'s are
distinct nonnegative integers. In particular, if $s=1$ then $n=2^w \ge 16$ and
$2^{[(n-s-1)/2]}=2^{(n-2)/2}>n-2$. This proves the property (b) if $s=1$.

If $s=2$ then
$$2^{[(n-s-1)/2]}=2^{(n-4)/2}=\frac{1}{2}2^{(n-2)/2}>(n-2)$$ 
if $n-2>8$. This proves the property (b) in the case of $n=2^{w_1}+2^{w_2}>10$. The remaining case $n=10$  follows from Tables in \cite{Atlas}.

Further we assume that $s \ge 3$. In particular, $n \ge 14$ and $n-2$ is {\sl not} a power of $2$. Since $n$ is even, all $w_i \ge 1$ and  $n \ge 2(2^s-1)$. If
$n = 2(2^s-1)$ then $n-2=2(2^s-2)$ is {\sl not} divisible by $2^3$. Therefore
if $[(n-s-1)/2]\ge 3$ then the property (b) holds. On the other hand, if
$[(n-s-1)/2]\le 3$ then $n-s-1 \le 2 \cdot 3+1$, i.e.,
$2(2^s-1)-s-1 \le 7$ which implies that $s<3$. We get a contradiction which, in turn, implies that $n>2(2^s-1)$ and therefore $n>2^{s+1}$.
This implies easily that $n \ge 2^{s+1}+6 \ge 20$. In particular,
$s<\log_2(n)-1$ and $2^{[(n-s-1)/2]} \ge 2^{[n-\log_2(n)]/2}>(2^n/n)^{1/2}/2$.

Since $n-2$ is not a power of $2$,  there are
no proper projective representations of $\An$ of dimension $2g=n-2$ in characteristic $0$ if 
$2^{[(n-s-1)/2]}>(n-2)/2.$
Clearly, this inequality holds if
$(2^n/n)^{1/2}/2>(n-2)/2$.
But this inequalty is equivalent to
$$2^n > n(n-2)^2$$
which holds for all $n \ge 20$.
\end{proof}


\begin{thebibliography}{99}
\bibitem{AbhMA} S. S. Abhyankar, {\em Alternating group coverings of the
affine line for characteristic greater than two}. Math. Ann. {\bf 296}
(1993), 63--68.

\bibitem{Atlas} J. H. Conway, R. T. Curtis, S. P. Norton, R. A. Parker, R.
A. Wilson, Atlas of finite groups. Clarendon Press, Oxford, 1985.

\bibitem{JK} G. James, A. Kerber, The representation theory of the symmetric
group. Addison Wesley Publishing Company, Reading, MA 1981.




\bibitem{Mori1} Sh. Mori, {\em The endomorphism rings of some abelian
varieties}. Japanese J. Math,  {\bf 2}(1976), 109--130.

\bibitem{Mori2} Sh. Mori, {\em The endomorphism rings of some abelian
varieties}. II, Japanese J. Math,  {\bf 3}(1977), 105--109.

\bibitem{MumfordAV} D. Mumford, Abelian varieties, Second edition,
 Oxford University Press, London, 1974.

\bibitem{Serre} J.-P. Serre, Topics in Galois Theory. Jones and Bartlett
Publishers, Boston-London, 1992.

\bibitem{Wagner} A. Wagner, {\em The faithful linear representations of
$\Sn$ and $\An$ over a field of characteristic} $2$. Math. Z. {\bf 151}
(1976), 127--137.

\bibitem{Wagner2} A. Wagner, {\em An observation on the degrees of the
symmetric and alternating group over an arbitrary field}. Arch. Math. 
{\bf 39} (1977), 583--589.

\bibitem{ZarhinMRL} Yu. G. Zarhin, {\em Hyperelliptic jacobians without
complex multiplication}. Math. Res. Letters {\bf 7}(2000), 123--132.

\bibitem{ZarhinTexel}Yu. G. Zarhin {\em Hyperelliptic jacobians and modular
representations}. In: Moduli of abelian varieties (C. Faber, G. van der
Geer, F. Oort, eds.), Birkh\"auser, to appear.

\end{thebibliography}
\end{document}